\documentstyle[leqno]{article}

\newtheorem{th}{Theorem}[section]
\newtheorem{prop}[th]{Proposition}

\newcounter{defin}[section]
\renewcommand{\thedefin}{\thesection.\arabic{defin}}
\newcounter{ex}[section]
\renewcommand{\theex}{\thesection.\arabic{ex}}
\newcounter{rem}[section]
\renewcommand{\therem}{\thesection.\arabic{rem}}
\date{}
\author{Mircea Neagu and Constantin Udri\c ste}
\title{Multi-Time Dependent Sprays and\\
Harmonic Maps on $J^1(T,M)$}
\begin{document}
\maketitle
\begin{abstract}
It is known that the jet fibre bundle of order one $J^1(T,M)$ is a basic
object in the study of classical and quantum field theories. In order to
develope a subsequent multi-time dependent Lagrangian theory of physical
fields on $J^1(T,M)$, we need to
generalize  the main geometrical objects used in  the classical rheonomic
Lagrangian theory.
In this direction, Section 1 presents the main properties of the differentiable
structure of the jet fibre bundle of order one. Section 2 studies
an important collection of geometrical objects on $J^1(T,M)$ as d-tensors,
temporal and spatial sprays and the harmonic maps induced by these sprays,
which naturally generalize analogous objects on $R\times TM$, the natural
house of the time-dependent Lagrangian field theory \cite{6}. Section 3
studies the nonlinear connections $\Gamma$ on $J^1(T,M)$, and discuss their
relation with the temporal and spatial sprays. Section 4 opens the problem
of prolongation of vector fields from $T\times M$ to 1-jet space
$J^1(T,M)$, using adapted bases.
\end{abstract}

\section {The jet fibre bundle $J^1(T,M)$}

\hspace{5mm} Let us consider the smooth manifolds $T$ and $M$ of
dimension $p$, respectively $n$, coordinated by $(t^\alpha)_{\alpha=\overline
{1,p}}$, respectively $(x^i)_{i=\overline{1,n}}$. We remark that, throughout
this paper, the set $\{1,2,\ldots,p\}$ is indexed by $\alpha,\beta,\gamma,
\ldots$, and the set $\{1,2,\ldots,n\}$ is indexed by $i,j,k,\ldots$.

Now, let $(t_0,x_0)$ be an arbitrary point of the product manifold $T\times M$.
We denote $C^\infty(T,M)$ the set of all smooth maps between $T$ and $M$
and define the equivalence relation
\begin{equation}
f\sim_{(t_{\scriptscriptstyle{0}},x_{\scriptscriptstyle{0}})}g
\Leftrightarrow\left\{
\begin{array}{l}
f(t_0)=g(t_0)=x_0\\
df_{t_{\scriptscriptstyle{0}}}=dg_{t_{\scriptscriptstyle{0}}},
\end{array}\right.
\end{equation}
on $C^\infty(T,M)$.
For every $f,g\in C^\infty(T,M)$, the relation
$f\sim_{(t_{\scriptscriptstyle{0}},x_{\scriptscriptstyle{0}})}g$ can be expressed locally by
\begin{equation}
\left\{\begin{array}{l}\medskip
x^i(t^\beta_0)=y^i(t^\beta_0)=x^i_0\\
\displaystyle{{\partial x^i\over\partial t^\alpha}(t^\beta_0)={\partial y^i\over\partial t^
\alpha}(t^\beta_0)},
\end{array}\right.
\end{equation}
where $t^\beta(t_0)=t^\beta_0,\;x^i(x_0)=x^i_0,\;x^i=x^i\circ f$ and $y^i=x^i
\circ g$. The equivalence class of a smooth map $f\in C^\infty(T,M)$ is
denoted by
$[f]_{(t_{\scriptscriptstyle{0}},x_{\scriptscriptstyle{0}})}=\{g\in C^\infty
(T,M)\;\vert\;g\sim_{(t_{\scriptscriptstyle{0}},x_{\scriptscriptstyle{0}})}
\nolinebreak f\}$.
If the quotient
$J^1_{t_{\scriptscriptstyle{0}},x_{\scriptscriptstyle{0}}}(T,M)=C^\infty(T,M)
\left/_{\sim_{(t_{\scriptscriptstyle{0}},x_{\scriptscriptstyle{0}})}}\right.$
is the factorization by the equivalence relation "$\sim_{(t_{\scriptscriptstyle{0}},
x_{\scriptscriptstyle{0}})}$", we build the total space of the 1-jet set, taking
\begin{equation}
J^1(T,M)=\bigcup_{(t_{\scriptscriptstyle{0}},x_{\scriptscriptstyle{0}})\in
T\times M}J^1_{t_{\scriptscriptstyle{0}},x_{\scriptscriptstyle{0}}}(T,M).
\end{equation}

Let us organize the total space of 1-jets $J^1(T,M)$ as a vector bundle over
the base space $T\times M$. To do this fact, we start with a smooth map
$f\in C^\infty(T,M)$,
$(t^1,\ldots,t^p)\to(x^1(t^1,\ldots,t^p),
\ldots,x^n(t^1,\ldots,t^p))$, and expand the maps $x^i$ using Taylor formula
around the point $(t^1_0,\ldots,t^p_0)\in R^p$. We obtain
$$
x^i(t^1,\ldots,t^p)=x^i_0+(t^\alpha-t^\alpha_0){\partial x^i\over\partial t^\alpha}
(t^1_0,\ldots,t^p_0)+{\cal O}(2),\quad\Vert(t^1-t^1_0,\ldots,t^p-t^p_0)\Vert<
\varepsilon.
$$
Considering the smooth map $\tilde f\in C^\infty(T,M)$ defined by the local
functions set
$$
\tilde x^i(t^1,\ldots,t^p)=x^i_0+(t^\alpha-t^\alpha_0){\partial x^i\over\partial t^\alpha}
(t^1_0,\ldots,t^p_0),\quad\Vert(t^1-t^1_0,\ldots,t^p-t^p_0)\Vert<\varepsilon,
$$
we deduce that $\tilde f\sim_{(t_{\scriptscriptstyle{0}},x_{\scriptscriptstyle
{0}})}f$, that is, the linear affine approximation $\tilde f$ of $f$ is a
convenient representative of the equivalence class
$[f]_{(t_{\scriptscriptstyle{0}},x_{\scriptscriptstyle{0}})}$.

Let $\pi:J^1(T,M)\to T\times M$ be the projection defined by
$\pi([f]_{(t_{\scriptscriptstyle{0}},x_{\scriptscriptstyle{0}})})=(t_0,f(t_0))$.
It is obvious that the map $\pi$ is well defined and surjective. Using this
projection, for every local chart $U\times V\subset T\times M$ on the
product manifold $T\times M$, we can define the bijection
$$
\Phi_{U\times V}:\pi^{-1}(U\times V)\to U\times V\times R^{np},
$$
setting $\Phi_{U\times V}([f]_{(t_{\scriptscriptstyle{0}},x_{\scriptscriptstyle
{0}})})=(t_0,x_0,\displaystyle{{\partial x^i\over\partial t^
\alpha}(t^\beta_0))},\;x_0=f(t_0).$

In conclusion, the 1-jet set $J^1(T,M)$ can be endowed with a differentiable
structure of dimension $p+n+pn$, such that the maps $\Phi_{U\times V}$ to be
diffeomorphisms. We emphasize that the local coordinates on $J^1(T,M)$ are
$(t^\alpha,x^i,x^i_\alpha)$, where
\begin{equation}\label{coor}
\left\{\begin{array}{l}
t^\alpha([f]_{(t_{\scriptscriptstyle{0}},x_{\scriptscriptstyle{0}})})=
t^\alpha(t_0)\\
x^i([f]_{(t_{\scriptscriptstyle{0}},x_{\scriptscriptstyle{0}})})=x^i(x_0)\\
x^i_\alpha([f]_{(t_{\scriptscriptstyle{0}},x_{\scriptscriptstyle{0}})})=
\displaystyle{{\partial x^i\over\partial t^\alpha}(t^\beta_0).}
\end{array}\right.
\end{equation}

In the above coordinates on $J^1(T,M)$, the projection $\pi:J^1(T\times M)\to T
\times M$ has the local expression $\pi(t^\alpha,x^i,x^i_\alpha)=(t^\alpha,x^i)$.
Moreover, the differential $\pi_*$ of the map $\pi$ is locally determined
by the Jacobi matrix
$$
\left(\begin{array}{ccc}
\delta_{\alpha\beta}&0&0\\
0&\delta_{ij}&0
\end{array}\right)\in M_{p+n,p+n+pn}.
$$
It  follows that $\pi_*$ is a surjection (rank $\pi_*=p+n$), and therefore
the projection $\pi$ is a submersion. Consequently, the 1-jet total space
$J^1(T,M)$ becomes a vector bundle over the base space $T\times M$, having
the fibre type $R^{pn}$.

Using \ref{coor}, by a simple direct calculation, we obtain
\begin{prop}
The local coordinate transformations
$(t^\alpha,x^i,x^i_\alpha)\longleftrightarrow(\tilde t^\alpha,\tilde x^i,\tilde
x^i_\alpha)$
of the 1-jet vector bundle $E=J^1(T,M)$ are given by
\begin{equation}
\left\{\begin{array}{l}
\tilde t^\alpha=\tilde t^\alpha(t^\beta)\\
\tilde x^i=\tilde x^i(x^j)\\
\tilde x^i_\alpha=\displaystyle{{\partial\tilde x^i\over\partial x^j}{\partial
t^\beta\over\partial\tilde t^\alpha}x^j_\beta,}
\end{array}\right.
\end{equation}
where $\det(\partial\tilde t^\alpha/\partial t^\beta)\neq 0$ and
$\det(\partial\tilde x^i/\partial x^j)\neq 0$. Consequently, $E$ is always an
orientable manifold.
\end{prop}

Let  us consider the canonical basis $\displaystyle{\left\{{\partial\over
\partial t^\alpha},{\partial\over\partial x^i},{\partial\over\partial x^i_\alpha}
\right\}}$ of vector fields on $E$ and $\{dt^\alpha,dx^i,dx^i_\alpha\}$
its dual basis of 1-forms.

\begin{prop}
Changing the coordinates on $E$, the  following transformation rules are true:
\begin{equation}
\left\{\begin{array}{l}\displaystyle{
{\partial\over\partial t^\alpha}={\partial\tilde t^\beta\over\partial t^\alpha}
{\partial\over\partial\tilde t^\beta}+{\partial\tilde x^j_\beta\over\partial t^\alpha}
{\partial\over\partial\tilde x^j_\beta}}\\\displaystyle{
{\partial\over\partial x^i}={\partial\tilde x^j\over\partial x^i}
{\partial\over\partial\tilde x^j}+{\partial\tilde x^j_\beta\over\partial x^i}
{\partial\over\partial\tilde x^j_\beta}}\\\displaystyle{
{\partial\over\partial x^i_\alpha}={\partial\tilde x^j\over\partial x^i}
{\partial t^\alpha\over\partial\tilde t^\beta}{\partial\over\partial\tilde
x^j_\beta},}
\end{array}\right.
\end{equation}
\begin{equation}
\left\{\begin{array}{l}\medskip\displaystyle{
dt^\alpha={\partial t^\alpha\over\partial\tilde t^\beta}d\tilde t^\beta}\\
\medskip\displaystyle{dx^i={\partial x^i\over\partial\tilde x^j}d\tilde x^j}\\
\displaystyle{
dx^i_\alpha={\partial x^i_\alpha\over\partial\tilde t^\beta}d\tilde t^\beta+
{\partial x^i_\alpha\over\partial\tilde x^j}d\tilde x^j +{\partial x^i\over
\partial\tilde x^j}{\partial\tilde t^\beta\over\partial t^\alpha}d\tilde x^j_
\beta.}
\end{array}\right.
\end{equation}
\end{prop}
{\bf \underline{Some physical apsects.}}\medskip

At the end of this Section, we should like to expose certain physical
aspects of the jet vector bundle of order one that we consider very eloquent
for the subsequent theory.

Thus, from physical point  of view, we regard the space $T$ as a "{\it temporal} "
manifold or a "{\it multi-time} " while the manifold $M$ is regarded as a
"{\it spatial} " one. The vector bundle $J^1(T,M)\to T\times M$ is
regarded as a {\it bundle of configurations}, in mechanics terms,  and its
elements $[f]$ are regarded as classes of "{\it parametrized sheets}".

In order to motivate the terminology used, we study more deeply the
jet vector bundle of order one, in the particular case $T=R$ (i. e., the usual time axis
represented by the set of real numbers).
Let us suppose that $J^1(R,M)\equiv R\times TM$ is coordinated by $(t,x^i,y^i)$.
The gauge group of the bundle
\begin{equation}
\pi:J^1(R,M)\to R\times M,\;\;(t,x^i,y^i)\to
(t,x^i),
\end{equation}
is given by
\begin{equation}\label{G_1}
\left\{\begin{array}{l}
\tilde t=\tilde t(t)\\
\tilde x^i=\tilde x^i(x^j)\\
\displaystyle{\tilde y^i={\partial\tilde x^i\over\partial x^j}{dt\over
d\tilde t}y^j.}
\end{array}\right.
\end{equation}
We  remark that the form of this gauge group stands out by the {\it relativistic}
character of the time $t$. For that reason, we consider that the jet fibre
bundle of order one $J^1(R,M)$ is the natural bundle  of configurations of the
{\it relativistic rheonomic Lagrangian mechanics} \cite{11}.

Comparatively, in the {\it classical rheonomic Lagrangian
mechanics} \cite{6}, the bundle of configurations is the fibre bundle
\begin{equation}
\pi: R\times TM\to M,\;(t,x^i,y^i)\to (x^i),
\end{equation}
whose geometrical invariance group is
\begin{equation}\label{G_2}
\left\{\begin{array}{l}
\tilde t=t\\
\tilde x^i=\tilde x^i(x^j)\\
\displaystyle{\tilde y^i={\partial\tilde x^i\over\partial x^j}y^j.}
\end{array}\right.
\end{equation}
Obviously, the structure of the gauge group \ref{G_2} emphasizes the  {\it
absolute}  character of the time $t$ from the classical rheonomic Lagrangian
mechanics. At the same time, we point out that the gauge group \ref{G_2} is
a subgroup of \ref{G_1}. In other words, the gauge group of the jet bundle of
order one, from the relativistic rheonomic Lagrangian mechanics, is more general
than that used in the classical rheonomic Lagrangian mechanics, which ignores
the temporal reparametrizations.

Finally, we invite the reader to compare both the  classical and relativistic
rheonomic  Lagrangian mechanics developed in \cite{6} and \cite{11}.

\section{d-Tensors. Multi-time dependent sprays. Harmonic maps}

\setcounter{equation}{0}
\hspace{5mm} It is well known the importance of the tensors in the development
of a geometry on a fibre bundle. In the study of the 1-jet fibre bundle,
a central role is played by the {\it distinguished tensors} or {\it d-tensors}.
\medskip\\
\addtocounter{defin}{1}
{\bf Definition \thedefin} A geometrical object
$D=(D^{\alpha i(j)(\nu)\ldots}_{\gamma k(\beta)(l)\ldots})$, on the 1-jet
vector bundle $E$, whose local components verify the following rules of
transformation,
\begin{equation}
D^{\alpha i(j)(\nu)\ldots}_{\gamma k(\beta)(l)\ldots}=\tilde
D^{\delta p(m)(\eta)\ldots}_{\varepsilon r(\mu)(s)\ldots}
{\partial t^\alpha\over\partial\tilde t^\delta}
{\partial x^i\over\partial\tilde x^p}
{\partial x^j\over\partial\tilde x^m}
{\partial\tilde t^\mu\over\partial t^\beta}
{\partial\tilde t^\varepsilon\over\partial t^\gamma}
{\partial\tilde x^r\over\partial x^k}
{\partial\tilde x^s\over\partial x^l}
{\partial t^\nu\over\partial\tilde t^\eta}\ldots,
\end{equation}
is called a {\it d-tensor field}.\medskip\\
\addtocounter{rem}{1}
{\bf Remarks \therem} i) The utilization of parentheses for certain indices of the local
components $D^{\alpha i(j)(\nu)\ldots}_{\gamma k(\beta)(l)\ldots}$ will be
motivated at the end of the Section 3 of this paper, before the introduction
of a nonlinear connection $\Gamma$ on $E$ together with its adapted bases of
vector and covector fields (see Remark 3.2).

ii) A d-tensor field $D$ on $E=J^1(T,M)$ can be viewed
like an object defined on $T\times M$ which depends on {\it partial derivatives}
or  {\it partial directions} $x^i_\alpha$.\medskip\\
\addtocounter{ex}{1}
{\bf Examples \theex} i) If $L:E\to R$ is a multi-time Lagrangian function with partial
derivatives of order one, the local components
\begin{equation}\label{mdt}
G^{(\alpha)(\beta)}_{(i)(j)}={1\over 2}{\partial^2L\over\partial x^i_\alpha
\partial x^j_\beta}
\end{equation}
represent a d-tensor field on $E$. We point out  that taking $T=R$ and $L$ a regular
time-dependent Lagrangian, the d-tensor field $G^{(1)(1)}_{(i)(j)}(t,x^i,y^i)$
is a natural generalization  of that so-called metrical
d-tensor field $g_{ij}(t,x,y)$ of  a classical  rheonomic Lagrange space
$RL^n=(M,L(t,x^i,y^i))$ \cite{6}.

ii) The geometrical object $\mbox{\bf C}=(\mbox{\bf C}^{(i)}_{(\alpha)})$,
where $\mbox{\bf C}^{(i)}_{(\alpha)}=x^i_\alpha$, represent a d-tensor field
on  $E$. This  is called the {\it canonical Liouville d-tensor} on the 1-jet
vector bundle $E$. We  emphasize that  this d-tensor field  naturally generalizes
the Liouville d-vector field {\bf C}$=y^i\displaystyle{{\partial\over\partial
y^i}}$ used in \cite{6}.

iii) Let $h_{\alpha\beta}$ be a semi-Riemannian metric on the temporal
manifold $T$. The geometrical object $L=(L^{(i)}_{(\alpha)\beta\gamma})$,
where $L^{(i)}_{(\alpha)\beta\gamma}=h_{\beta\gamma}x^i_\alpha$, is a d-tensor
field which is called the {\it Liouville d-tensor associated to the metric $h$}.

iv) Using the preceding metric $h$, we construct the d-tensor
$J=(J^{(i)}_{(\alpha)\beta j})$, where $J^{(i)}_{(\alpha)\beta j}=h_{\alpha\beta}
\delta^i_j$. This d-tensor is called the {\it $h$-normalization d-tensor}
of the jet bundle $E$. Note that the $h$-normalization d-tensor of $J^1(T,M)$
is a natural generalization of the tangent structure $J$ from the Lagrange
geometry \cite{6}.

It is obvious that any d-tensor on $E$ is a tensor on $E$. Conversely, this
is not true. As examples, we will build two tensors which are not d-tensors.
We refer to notions of temporal  and spatial sprays which allow the generalization
of the notion of time-dependent spray used  in \cite{6},  \cite{14}.\medskip\\
\addtocounter{defin}{1}
{\bf Definition \thedefin} A global tensor $H$, expressed locally by
\begin{equation}
H=\delta^\beta_\alpha dt^\alpha\otimes{\partial\over\partial t^\beta}-
2H^{(j)}_{(\beta)\alpha}dt^\alpha\otimes{\partial\over\partial x^j_\beta},
\end{equation}
is called a {\it temporal spray} on $E$.\medskip

Taking into account  that  a temporal spray is a global tensor on $E$, by a
direct  calculation, we deduce
\begin{prop}
i) The components $H^{(j)}_{(\beta)\alpha}$ of the temporal spray $H$
transform by the rules
\begin{equation}\label{tts}
2\tilde H^{(k)}_{(\mu)\gamma}=2H^{(j)}_{(\beta)\alpha}{\partial t^\alpha\over
\partial\tilde t^\gamma}{\partial\tilde x^k\over\partial x^j}
{\partial t^\beta\over\partial\tilde t^\mu}-{\partial t^\alpha\over\partial
\tilde t^\gamma}{\partial\tilde x^k_\mu\over\partial t^\alpha}.
\end{equation}

ii) Conversely, to give a temporal spray on $E$ is equivalent to give a set
of local functions $H=(H^{(j)}_{(\beta)\alpha})$ which transform by \ref{tts}.

iii) The global tensor
$$
H=H^\beta_\alpha dt^\alpha\otimes{\partial\over\partial t^\beta}-
2H^{(j)}_{(\beta)\alpha}dt^\alpha\otimes{\partial\over\partial x^j_\beta}
$$
is a temporal spray iff
$J^{(j)}_{(\beta)\alpha i}H^\alpha_\gamma=J^{(j)}_{(\beta)\gamma i}$,
where $J$ is the normalization d-tensor of the fibre bundle $E$ associated
to an arbitrary semi-Riemannian metric $h$ on $T$.
\end{prop}

The previous  proposition allows us to offer the following important example
of temporal spray. The importance of this temporal spray is determined  by
its using in the  description of the classical harmonic maps between two
semi-Riemannian manifolds \cite{3}.\medskip\\
\addtocounter{ex}{1}
{\bf Example \theex} Using the  transformation rules of the Christoffel
symbols $H^\alpha_{\beta\gamma}$ attached to a  semi-Riemannian metric
$h_{\alpha\beta}$ on $T$, we deduce that the components $2H^{(j)}_{(\beta)
\alpha}=-H_{\alpha\beta}^\gamma x^j_\gamma$ represent a temporal spray on $E$.
This  is called the {\it canonical temporal spray associated to the metric
$h$.}\medskip\\
\addtocounter{defin}{1}
{\bf Definition \thedefin}  A global tensor $G$, locally defined by
\begin{equation}
G=x^i_\alpha dt^\alpha\otimes{\partial\over\partial x^i}-
2G^{(j)}_{(\beta)\alpha}dt^\alpha\otimes{\partial\over\partial x^j_\beta},
\end{equation}
is called a {\it spatial spray} on $E$.

As in the case of the temporal spray, we can prove without difficulties
the following statements.
\begin{prop}
The components $G^{(j)}_{(\beta)\alpha}$ of the spatial spray $G$ transform
by the rules
\begin{equation}\label{tss}
2\tilde G^{(k)}_{(\mu)\gamma}=2G^{(j)}_{(\beta)\alpha}{\partial t^\alpha\over
\partial\tilde t^\gamma}{\partial\tilde x^k\over\partial x^j}
{\partial t^\beta\over\partial\tilde t^\mu}-{\partial x^i\over\partial
\tilde x^j}{\partial\tilde x^k_\mu\over\partial x^i}\tilde x^j_\gamma.
\end{equation}

ii) To give a spatial spray is equivalent to give a
set of local functions $G=(G^{(j)}_{(\beta)\alpha})$ which change by the
law \ref{tss}.

iii) A global tensor on $E$, defined locally by
$$
G=G^i_\alpha dt^\alpha\otimes{\partial\over\partial x^i}-
2G^{(j)}_{(\beta)\alpha}dt^\alpha\otimes{\partial\over\partial x^j_\beta},
$$
is a spatial spray iff $J^{(j)}_{(\beta)\alpha i}G^i_\gamma=L^{(j)}_{(
\beta)\alpha\gamma}$, where $J$ (resp. $L$) is the normalization (resp.
Liouville) d-tensor associated to an arbitrary semi-Riemannian metric $h$.
\end{prop}
\addtocounter{ex}{1}
{\bf Example \theex} If $\gamma^i_{jk}$ are the Christoffel symbols of a
semi-Riemannian metric $\varphi_{ij}$ on the spatial manifold $M$, the local
coefficients $2G^{(j)}_{(\beta)\alpha}=\gamma^j_{kl}x^k_\alpha x^l_\beta$ define a spatial spray which
is called the {\it canonical spatial spray associated to the metric $\varphi$.}
We point out that this spatial  spray is also used  in the description of  the
classical harmonic  maps between  two semi-Riemannian manifolds \cite{3}.
\medskip\\
\addtocounter{defin}{1}
{\bf Definition \thedefin} A pair $(H,G)$, which consists of a temporal spray
and a spatial one, is called  a {\it  multi-time dependent spray} on $E$.\medskip

To characterize the multi-time dependent sprays on $E$ and to underline again
the importance of the canonical temporal and spatial sprays attached to
the metrics $h$ and $\varphi$, we prove  the following
\begin{th}\label{char}
Let $(T,h),\;(M,\varphi)$ be semi-Riemannian manifolds
and let $H=\nolinebreak(H^{(i)}_{(\alpha)\beta})$ (resp. $G=(G^{(i)}_{(\alpha)\beta})$)
be an arbitrary temporal (resp. spatial) spray on $E$. In these
conditions, we have
\begin{equation}
\left\{\begin{array}{l}\medskip
H^{(i)}_{(\alpha)\beta}=-{1\over 2}H^\gamma_{\alpha\beta}x^i_\gamma+
D^{(i)}_{(\alpha)\beta}\\
G^{(i)}_{(\alpha)\beta}={1\over 2}\gamma^i_{jk}x^j_\alpha x^k_\beta+
F^{(i)}_{(\alpha)\beta},
\end{array}\right.
\end{equation}
where $D^{(i)}_{(\alpha)\beta},\;F^{(i)}_{(\alpha)\beta}$ are certain d-tensors
on $E$.
\end{th}
{\bf Proof.} The theorem comes from the following true statements:

i) An affine combination of temporal (spatial) sprays is a temporal (spatial)
spray,

ii) The product between a scalar and a temporal (spatial) spray is a temporal
(spatial) spray,

iii) The difference between two temporal (spatial) sprays is a d-tensor.
\rule{5pt}{5pt}\medskip

In order to  generalize the  notion of  path of a spray  from Lagrangian geometry,
we fixe $h_{\alpha\beta}$ a semi-Riemannian
metric on the temporal manifold $T$. In this context, we give the following
\medskip\\
\addtocounter{defin}{1}
{\bf Definition \thedefin} A geometrical object $H=(H^k)$ (resp. $G=(G^k)$)
is called a {\it temporal} (resp. {\it spatial} ) {\it $h$-spray} if the local
components modify by the rules
\begin{equation}
2\tilde H^k=2H^j{\partial\tilde x^k\over\partial x^j}-\tilde h^{\gamma\mu}
{\partial t^\alpha\over\partial\tilde t^\gamma}
{\partial\tilde x^k_\mu\over\partial t^\alpha},
\end{equation}
repectively
\begin{equation}\label{thss}
2\tilde G^k=2G^j{\partial\tilde x^k\over\partial x^j}-\tilde h^{\gamma\mu}
{\partial x^i\over\partial\tilde x^j}
{\partial\tilde x^k_\mu\over\partial x^i}\tilde x^j_\gamma.
\end{equation}\medskip\\
\addtocounter{ex}{1}
{\bf Example \theex} Starting with $H=(H^{(i)}_{(\alpha)\beta})$ (resp. $G=(G^{(i)}
_{(\alpha)\beta})$) like a temporal (resp. spatial ) spray, the entity $H=(H^i)$ (resp.
$G=(G^i)$), where $H^i=h^{\alpha\beta}H^{(i)}_{(\alpha)\beta}$ (resp.
$G^i=h^{\alpha\beta}G^{(i)}_{(\alpha)\beta}$), represents a temporal (resp.
spatial ) $h$-spray which will be called the {\it $h$-trace of the temporal}
(resp. {\it spatial} ) spray $H$ (resp. $G$). Particularly, the components
$H^k=-h^{\alpha\beta}H^\gamma_{\alpha\beta}x^k_\gamma$ (resp. $G^k=h^{\alpha
\beta}\gamma^k_{ij}x^i_\alpha x^j_\beta$) represent the {\it canonical temporal}
(resp. {\it spatial} ) {\it h-spray attached to the metric } $h$ (resp.
$\varphi$).\medskip

The previous example show that the $h$-trace of a temporal or a spatial spray
represents a temporal or a spatial $h$-spray. Conversely, we prove the following
\begin{th}
If $\dim T=1$, any temporal (spatial) $h$-spray is the
$h$-trace of a unique temporal (spatial ) spray.
\end{th}
{\bf Proof.} Let $G=(G^k)$ be a spatial $h$-spray. We denote $G^{(k)}_{(1)1}=
h_{11}G^k$. Obviously, the relation $G^k=h^{11}G^{(k)}_{(1)1}$ is true. In these
conditions, using the transformation rules \ref{thss}, we deduce
$$
2\tilde G^{(k)}_{(1)1}=2G^{(j)}_{(1)1}{\partial\tilde x^k\over\partial x^j}
\left({d\tilde t\over dt}\right)^2-{dt\over d\tilde t}{dy^k\over dt}.
$$
This means that $G=(G^{(k)}_{(1)1})$ is a spatial spray. The uniqueness is clear.

By analogy, we treat the case of the temporal $h$-sprays, taking
$H^k=h^{11}H^{(k)}_{(1)1}$, where $H^{(k)}_{(1)1}=h_{11}H^k$. \rule{5pt}{5pt}
\medskip\\
\addtocounter{rem}{1}
{\bf Remark \therem} The previous theorem shows that, in the case
$\dim T=1$, there is a 1-1 corespondence between sprays and $h$-sprays while,
for $\dim T\geq 2$, this statement is not true.\medskip

In the sequel, let us fixe a temporal  spray  $H=(H^{(i)}_{(\alpha)\beta})$ and
\linebreak a spatial  spray  $G=(G^{(i)}_{(\alpha)\beta})$ on $E$.
The following notions show that the 1-jet fibre bundle is the natural house
for important objects with geometrical and physical meaning.\medskip\\
\addtocounter{defin}{1}
{\bf Definition \thedefin} A solution $f\in C^\infty(T,M)$ of the PDEs system  of order two
\begin{equation}\label{afmap}
x^i_{\alpha\beta}+G^{(i)}_{(\alpha)\beta}+G^{(i)}_{(\beta)\alpha}+
H^{(i)}_{(\alpha)\beta}+H^{(i)}_{(\beta)\alpha}=0,
\end{equation}
where the map $f$ is locally expressed by $(t^\alpha)\to (x^i(t^\alpha))$ and
$\displaystyle{x^i_{\alpha\beta}={\partial^2x^i\over\partial t^\alpha\partial
t^\beta}}$, is called {\it an affine map of the multi-time dependent spray
$(H,G)$.}\medskip\\
\addtocounter{rem}{1}
{\bf Remark \therem} A reason which offers the naturalness of the notion of an
affine map of a multi-time dependent spray on $E$, is that, in the particular
case $T=R$, the equations of the affine maps generalize the equations of the
paths of a time-dependent spray from the rheonomic Lagrangian geometry \cite{2}.\medskip\\
\addtocounter{ex}{1}
{\bf Example \theex} Considering the canonical multi-time dependent spray
\begin{equation}
\left\{\begin{array}{l}\medskip
H^{(i)}_{(\alpha)\beta}=-{1\over 2}H^\gamma_{\alpha\beta}x^i_\gamma\\
G^{(i)}_{(\alpha)\beta}={1\over 2}\gamma^i_{jk}x^j_\alpha x^k_\beta,
\end{array}\right.
\end{equation}
the equations of the affine maps of this spray reduce to
\begin{equation}
x^i_{\alpha\beta}-H^\gamma_{\alpha\beta}x^i_\gamma+\gamma^i_{jk}x^j_\alpha
x^k_\beta=0,
\end{equation}
that is, the equations whose solutions are exactly the maps
$f\in C^\infty(T,M)$ which carry the geodesics of $(T,h_{\alpha\beta})$ into the
geodesics of the space $(M,\varphi_{ij})$.\medskip

Taking $h=(h_{\alpha\beta}(t))$ a temporal semi-Riemannian metric and doing
a contraction by $h^{\alpha\beta}$ in \ref{afmap}, we can introduce the next
\medskip\\
\addtocounter{defin}{1}
{\bf Definition \thedefin} A map $f\in C^\infty(T,M)$ is called a {\it harmonic
map of the  multi-time dependent spray $(H,G)$, with respect to the semi-Riemannian
metric $h$}, if $f$ is a solution of the PDEs system of order two
\begin{equation}
h^{\alpha\beta}\{x^i_{\alpha\beta}+2G^{(i)}_{(\alpha)\beta}+2H^{(i)}_{(\alpha)
\beta}\}=0.
\end{equation}\medskip\\
\addtocounter{ex}{1}
{\bf Example \theex} Particularly, in the case of the canonical multi-time dependent
spray of preceding example, we recover the classical notion of harmonic map
between the semi-Riemannian manifolds $(T,h)$ and $(M,\varphi)$ \cite{3}. This
fact points out that our generalization of the classical notion of harmonic
map is a natural one.\medskip\\
\addtocounter{rem}{1}
{\bf Remarks \therem} i) It is obvious that the affine map of a multi-time dependent
spray $(H,G)$ is a harmonic map of the same spray with respect to any
semi-Riemannian metric $h_{\alpha\beta}$ on the temporal space $T$.

ii) In the particular case $(T,h)=(R,\delta)$, the notions of harmonic map and
affine map identify. Consequently, both notions naturally generalize that
so-called a path of a time-dependent spray, used  in \cite{6}.\medskip

Let us denote $S^{(i)}_{(\alpha)\beta}=G^{(i)}_{(\alpha)\beta}+H^{(i)}_{
(\alpha)\beta}+{1\over 2}H^\gamma_{\alpha\beta}x^i_\gamma$ and $S^i=h^{\alpha
\beta}S^{(i)}_{(\alpha)\beta}$. Using the theorem \ref{char},  we deduce
that $S=(S^i)$ is a spatial $h$-spray. In this context, we obtain the
without difficulties the  following
\begin{th}
The equations of the harmonic maps of the multi-time  dependent spray $(G,H)$,
with respect to the semi-Riemannian metric $h$, can be rewritten in the Poisson
form
\begin{equation}\label{pois}
\Delta_hx^i+2S^i=0,
\end{equation}
where $\Delta_hx^i=h^{\alpha\beta}(x^i_{\alpha\beta}-H^\gamma_{\alpha\beta}x^i_
\gamma)$.
\end{th}
\addtocounter{rem}{1}
{\bf Remarks \therem} i) This theorem will play a  central role in the  development
of the  subsequent (generalized) metrical  multi-time Lagrange  theory of physical fields. In this sense, we will prove over there that the Euler-Lagrange equations
of a  multi-time dependent Lagrangian ${\cal L}=L\sqrt{\vert h\vert}$, where
$L:J^1(T,M)\to R$ is a {\it Kronecker  $h$-regular} Lagrange function,
can be written in the Poisson form \ref{pois}. Hence, the extremals of ${\cal L}$ can be  regarded
as {\it harmonic maps}, in our  sense, offering them a profound geometrical
and  physical character. For more details, see \cite{8}, \cite{10}.

ii) On the other  hand, the same theorem will be used to offer
a beautiful geometrical interpretation of solutions of PDEs, in metrical
multi-time Lagrangian geometry terms \cite{18}. In this fashion, we will
offer a final answer to the Udri\c ste-Neagu open problem \cite{9}, \cite{17},
whose essential physical aspects are presented in \cite{15}, \cite{16}.

\section{Nonlinear connections}

\setcounter{equation}{0}
\hspace{5mm} The form of the coordinate changes on $E=J^1(T,M)$ determines
complicated rules of transformation of the local components of diverse
geometrical objects of this space. This motivates the introduction
of {\it nonlinear connection}  which induces {\it adapted bases}.
These bases have the quality to simplify the transformation rules of the
components of the geometrical objects taken in study.

With a view to doing this, we take $u\in E$ and consider the differential
map\linebreak $\pi_{*,u}:T_uE\to T_{(t,x)}(T\times M)$  of the canonical
projection $\pi:E\to T\times M,\;\pi(u)=(t,x)$. At the
same time, let us consider the vector subspace $V_u=Ker\;\pi_{*,u}\subset T_uE$.
Because the map $\pi_{*,u}$ is a surjection, we have $\dim_RV_u=pn,\;\forall
u\in E$. Moreover, a basis in $V_u$ is determined by $\{{\partial\over\partial
x^i_\alpha}\}$. In conclusion, the map
\begin{equation}
{\cal V}:\;u\in E\to V_u\subset T_uE
\end{equation}
is a differential distribution which is called {\it the vertical distribution}
of the 1-jet fibre bundle $E$.\medskip\\
\addtocounter{defin}{1}
{\bf Definition \thedefin} A {\it nonlinear connection} on $E$ is a differential
distribution
\begin{equation}
{\cal H}:\;u\in E\to H_u\subset T_uE
\end{equation}
which verifies the relation
\begin{equation}
T_uE=H_u\oplus V_u,\;\forall u\in E.
\end{equation}
The distribution ${\cal H}$ is called the {\it horizontal distribution} on $E$.
\medskip\\
\addtocounter{rem}{1}
{\bf Remarks \therem} i) The above definition implies that $\dim_RH_u=p+n,\;\forall
u\in E$.

ii) The vector fields set ${\cal X}(E)$ can be decompose in the following direct
sum ${\cal X}(E)=\nolinebreak\Gamma({\cal H})\oplus\Gamma({\cal V})$, where $\Gamma({\cal
H})$ (resp. $\Gamma({\cal V})$) is the set of the sections on ${\cal H}$ (resp.
${\cal V}$).\medskip

Now, supposing that there is a nonlinear connection ${\cal H}$ on $E$, we have
the isomorphism
\begin{equation}
\pi_{*,u}\vert_{H_u}:\;H_u \to T_{\pi(u)}(T\times M),
\end{equation}
which allows us to prove the following
\begin{th}
i) There exist the unique horizontal vector fields
$\displaystyle{{\delta\over\delta t^\alpha},\;{\delta\over\delta x^i}\in
\Gamma({\cal H})}$, linearly independent, having the property
\begin{equation}
\pi_*\left({\delta\over\delta t^\alpha}\right)={\partial\over\partial t^\alpha},
\;\pi_*\left({\delta\over\delta x^i}\right)={\partial\over\partial x^i}.
\end{equation}

ii) The vector fields $\displaystyle{{\delta\over\delta t^\alpha}}$ and
$\displaystyle{{\delta\over\delta x^i}}$ can be uniquely written in the form
\begin{equation}
\left\{\begin{array}{l}\displaystyle{
{\delta\over\delta t^\alpha}={\partial\over\partial t^\alpha}-M^{(j)}_{(\beta)
\alpha}{\partial\over\partial x^j_\beta}}\\
\displaystyle{{\delta\over\delta x^i}={\partial\over\partial x^i}-N^{(j)}_{(\beta)
i}{\partial\over\partial x^j_\beta}}.
\end{array}\right.
\end{equation}

iii) The coefficients $M^{(j)}_{(\beta)\alpha}$ and $N^{(j)}_{(\beta)i}$
modify by the rules
\begin{equation}\label{lnlc}
\left\{\begin{array}{l}\medskip\displaystyle{
\tilde M^{(j)}_{(\beta)\mu}{\partial\tilde t^\mu\over\partial t^\alpha}=
M^{(k)}_{(\gamma)\alpha}{\partial \tilde x^j\over\partial x^k}
{\partial t^\gamma\over\partial\tilde t^\beta}-{\partial\tilde x^j_\beta\over
\partial t^\alpha}}\\
\displaystyle{\tilde N^{(j)}_{(\beta)k}{\partial\tilde x^k\over\partial x^i}=
N^{(k)}_{(\gamma)i}{\partial \tilde x^j\over\partial x^k}
{\partial t^\gamma\over\partial\tilde t^\beta}-{\partial\tilde x^j_\beta\over
\partial x^i}}.
\end{array}\right.
\end{equation}

iv) To give a nonlinear connection ${\cal H}$ on $E$ is equivalent to give
a set of local functions $\Gamma=(M^{(j)}_{(\beta)\alpha},N^{(j)}_
{(\beta)i})$ which transform by \ref{lnlc}.
\end{th}
\addtocounter{ex}{1}
{\bf Example \theex} Studying the transformation rules of the local components
\begin{equation}
\left\{\begin{array}{l}\medskip
M^{(j)}_{(\beta)\alpha}=-H^\gamma_{\alpha\beta}x^j_\gamma\\
N^{(j)}_{(\beta)i}=\gamma^j_{ik}x^k_\beta,
\end{array}\right.
\end{equation}
we conclude that $\Gamma_0=(M^{(j)}_{(\beta)\alpha},N^{(j)}_{(\beta)i})$
represents a nonlinear connection on $E$, which is called the {\it canonical
nonlinear connection attached to the semi-Riemannian metrics $h_{\alpha\beta}$
and $\varphi_{ij}$}.\medskip

Let us consider the 1-form $\delta x^i_\alpha=dx^i_\alpha+M^{(i)}_{(\alpha)
\beta}dt^\beta+N^{(i)}_{(\alpha)j}dx^j$. One easily deduces that the set of
1-forms $\{dt^\alpha, dx^i,\delta x^i_\alpha\}$ is a basis in the set of
1-forms.\medskip\\
\addtocounter{defin}{1}
{\bf Definition \thedefin} The basis $\displaystyle{\left\{{\delta\over\delta t^\alpha},
{\delta\over\delta x^i},{\partial\over\partial x^i_\alpha}\right\}\subset {\cal X}
(E)}$ and its  dual basis $\{dt^\alpha, dx^i,\delta x^i_\alpha\}\subset
{\cal X}^*(E)$ are called the {\it adapted bases} on $E$, determined by the
nonlinear connection $\Gamma$.\medskip

The big advantage of the adapted bases is that the transformation laws of its
elements are simple and natural.
\begin{prop}
The transformation laws of the elements of the adapted bases attached to the
nonlinear connection $\Gamma$ are
\begin{equation}\label{vab}
\left\{\begin{array}{l}\medskip
\displaystyle{{\delta\over\delta t^\alpha}={\partial\tilde t^\beta\over\partial t^\alpha}
{\delta\over\delta\tilde t^\beta}}\medskip\\
\displaystyle{{\delta\over\delta x^i}={\partial\tilde x^j\over\partial x^i}
{\delta\over\delta\tilde x^j}}\medskip\\
\displaystyle{{\partial\over\partial x^i_\alpha}={\partial\tilde x^j\over\partial x^i}
{\partial t^\alpha\over\partial\tilde t^\beta}{\delta\over\delta\tilde x^j_
\beta},}
\end{array}\right.
\end{equation}
\begin{equation}\label{cvab}
\left\{\begin{array}{l}\medskip
\displaystyle{dt^\alpha={\partial t^\alpha\over\partial\tilde t^\beta}
d\tilde t^\beta}\medskip\\
\displaystyle{dx^i={\partial x^i\over\partial\tilde x^j}d\tilde x^j}
\medskip\\
\displaystyle{\delta x^i_\alpha={\partial x^i\over\partial\tilde x^j}
{\partial\tilde t^\beta\over\partial t^\alpha}\delta\tilde x^j_\beta.}
\end{array}\right.
\end{equation}
\end{prop}
\addtocounter{rem}{1}
{\bf Remark \therem} The simple transformation rules \ref{vab} and \ref{cvab}
determine us to describe the objects with geometrical and physical meaning
from the subsequent (generalized) metrical multi-time Lagrange theory of physical
fields \cite{8}, \cite{10}, in adapted components. In a such prospect, we emphasize
that, using adapted bases of a  nonlinear connection $\Gamma$, a d-tensor
$D=(D^{\alpha i(j)(\mu)\ldots}_{\gamma k(\beta)(l)\ldots})$
on $E$ can be regarded as a global geometrical object, locally defined by
\begin{equation}
D=D^{\alpha i(j)(\mu)\ldots}_{\gamma k(\beta)(l)\ldots}{\delta\over\delta t^\alpha}
\otimes{\delta\over\delta x^i}\otimes{\partial\over\partial x^j_\beta}\otimes
dt^\gamma\otimes dx^k\otimes\delta x^l_\mu\otimes\ldots\;.
\end{equation}
The utilization of certain indices between parenthesis in the description of
the local components of the d-tensor $D$ is suitable for contractions. To
illustrate this fact, we consider, for example, the local components  of the
metrical d-tensor \ref{mdt} from  the example 2.1. These define the
geometrical object
\begin{equation}
G=G^{(\alpha)(\beta)}_{(i)(j)}\delta x^i_\alpha\otimes\delta x^j_\beta.
\end{equation}
On the other hand, considering the local components of the $h$-normalization
d-tensor $J^{(i)}_{(\alpha)\beta j}$, we obtain the representative object
\begin{equation}
\displaystyle{J=J^{(i)}_{(\alpha)\beta j}{\delta\over\delta x^i_\alpha}\otimes dt^\beta
\otimes dx^j}.
\end{equation}

Finally, let us study the relation between the notion of spray and the nonlinear
connection. In this context, the coefficients $M^{(j)}_{(\beta)\alpha}$ (resp.
$N^{(j)}_{(\beta)i}$) of the nonlinear connection $\Gamma$ are called
the {\it temporal} (resp. {\it spatial}) {\it nonlinear connection.} In this
terminology, using the transformation formulas \ref{tts}, \ref{tss} and
\ref{lnlc}, we can easily prove the following statements.
\begin{th}
i) If $M^{(i)}_{(\alpha)\beta}$ are the components
of a temporal nonlinear connection, then the components
\begin{equation}
\displaystyle {H^{(i)}_{(\alpha)\beta}={1\over2}M^{(i)}_{(\alpha)\beta}}
\end{equation}
represent a temporal spray.\medskip

ii) Conversely, if $H^{(i)}_{(\alpha)\beta}$ are the components of a temporal
spray, then
\begin{equation}
M^{(i)}_{(\alpha)\beta}=\nolinebreak 2H^{(i)}_{(\alpha)\beta}
\end{equation}
are the components of a temporal nonlinear connection.
\end{th}

\begin{th}
i) If $G^{(i)}_{(\alpha)\beta}$ are the components of
a spatial spray and $G^i=\nolinebreak h^{\alpha\beta}G^{(i)}_{(\alpha)\beta}$ represent the
$h$-trace of this spray, then the coefficients
\begin{equation}
N^{(i)}_{(\alpha)j}={\partial G^i\over\partial x^j_\gamma}h_{\gamma\alpha}
\end{equation}
represent a spatial nonlinear connection.

ii) Conversely, the spatial nonlinear connection $N^{(i)}_{(\alpha)j}$ induces
the spatial spray
\begin{equation}
2G^{(i)}_{(\alpha)\beta}=N^{(i)}_{(\alpha)j}x^j_\beta.
\end{equation}
\end{th}
\addtocounter{rem}{1}
{\bf Remark \therem} The previous theorems allow us to conclude that a
multi-time dependent spray $(H,G)$  induces
naturally a nonlinear connection $\Gamma$ on $E$, which is called the {\it canonical
nonlinear connection associated to the multi-time dependent spray $(H,G)$.}
We point out that the canonical nonlinear  connection $\Gamma$ attached  to the
multi-time dependent  spray $(H,G)$ is a natural generalization of
the canonical nonlinear connection $N$ induced by a time-dependent spray $G$ from the
classical rheonomic Lagrangian geometry \cite{6}.

\section{Jet prolongation of vector fields}

\setcounter{equation}{0}
\hspace{5mm} A general vector field $X^*$ on $J^1(T,M)$ can be written under
the form
$$
X^*=X^\alpha{\partial\over\partial t^\alpha}+X^i{\partial\over\partial x^i}
+X^{(i)}_{(\alpha)}{\partial\over\partial x^i_\alpha},
$$
where the components $X^\alpha,\;X^i\;X^{(i)}_{(\alpha)}$ are functions of $(t^\alpha,x^i,x^i_\alpha)$.

The prolongation of a vector field $X$ on $T\times M$ to a vector field on the
1-jet bundle $J^1(T,M)$ was solved by Olver \cite{12} in the following sense.
\medskip\\
\addtocounter{defin}{1}
{\bf Definition \thedefin} Let $X$ be a vector field on $T\times M$ with corresponding
(local) one-parameter group $\exp(\varepsilon X)$. The {\it 1-th prolongation}
of $X$, denoted by $pr^{(1)}X$, will be a vector field on the 1-jet space
$J^1(T,M)$, and is defined to be the infinitesimal generator of the corresponding
prolonged one-parameter group $pr^{(1)}[\exp(\varepsilon X)]$, i. e. ,
\begin{equation}
\displaystyle{[pr^{(1)}X](t^\alpha,x^i,x^i_\alpha)=\left.{d\over d\varepsilon}
\right\vert_{\varepsilon=0}pr^{(1)}[\exp(\varepsilon X)](t^\alpha,x^i,x^i_
\alpha).}
\end{equation}

In order to write the components of the prolongation, Olver used the {\it
$\alpha$-th total derivative} $D_\alpha$ of an arbitrary
function $f(t^\alpha,x^i)$ on $T\times M$, which is defined by the relation
\begin{equation}
\displaystyle{D_\alpha f={\partial f\over\partial t^\alpha}+{\partial f\over
\partial x^i}x^i_\alpha.}
\end{equation}

Thus, starting with $\displaystyle{X=X^\alpha(t,x){\partial\over\partial t^\alpha}+
X^i(t,x){\partial\over\partial x^i}}$ like a vector field on $T\times M$,
Olver introduced the 1-th prolongation of $X$ as the vector field
\begin{equation}
pr^{(1)}X=X+X^{(i)}_{(\alpha)}(t^\beta,x^j,x^j_\beta){\partial
\over\partial x^i_\alpha},
\end{equation}
where
$$
X^{(i)}_{(\alpha)}=D_\alpha X^i-(D_\alpha X^\beta)x^i_\beta={\partial X^i\over
\partial t^\alpha}+{\partial X^i\over\partial x^j}x^j_\alpha-\left({\partial
X^\beta\over\partial t^\alpha}+{\partial X^\beta\over\partial x^j}x^j_\alpha
\right)x^i_\beta.
$$

If we assume that is given a nonlinear connection $\Gamma=(M^{(i)}_{(\alpha)\beta},
N^{(i)}_{(\alpha)j})$ on $J^1(T,M)$, then the $\alpha$-th total
derivative used by Olver can be written as
\begin{equation}
D_\alpha f={\delta f\over\delta t^\alpha}+{\delta f\over\delta x^i}x^i_\alpha,
\end{equation}
and, consequently, $D_\alpha f$ represent the local components
of a distinguished 1-form on $J^1(T,M)$, which is expressed by $Df=(D_\alpha
f)dt^\alpha$.

Now, let there be given a vector field $X$ on $T\times M$. From a geometrical
point of view, we can define a 1-jet prolongation of $X$ as the {\it horizontal
lift} $X^H$ of $X$. This is defined by
\begin{equation}
X^H=X^\alpha{\delta\over\delta t^\alpha}+
X^i{\delta\over\delta x^i}=X-(M^{(j)}_{(\beta)\alpha}X^\alpha+N^{(j)}_{(\beta
)i}X^i){\partial\over\partial x^j_\beta}.
\end{equation}\medskip\\
{\bf \underline{Open problem}.}

Study the prolongations of vectors, 1-forms, tensors, $G$-structures from
$T\times M$ to $J^1(T,M)$.\medskip\\
{\bf Acknowledgments.} A version  of this paper was presented at Third
Conference of Balkan Society of Geometers, Politehnica University of
Bucharest, Romania, July 31-August 3, 2000. It is a pleasure for us to thank
to Prof. Dr. D. Opri\c s and to the reviewers of Journal of the London
Mathematical Society for their valuable comments upon the previous version of
this paper.

\begin{center}
University POLITEHNICA of Bucharest\\
Department of Mathematics I\\
Splaiul Independentei 313\\
77206 Bucharest, Romania\\
e-mail: mircea@mathem.pub.ro\\
e-mail: udriste@mathem.pub.ro
\end{center}


\begin{thebibliography}{20}
\bibitem{1}
G. S. Asanov, {\it Gauge-Covariant Stationary Curves on Finslerian and Jet
Fibrations and Gauge Extension of Lorentz Force}, Tensor N. S. , Vol 50
(1991), 122-137.
\bibitem{2}
L. A. Cordero, C. T. J. Dodson, M. de L\'eon, {\it Differential Geometry of
Frame Bundles}, Kluwer Academic Publishers, 1989.
\bibitem{3}
J. Eells, L. Lemaire, {\it A report on harmonic maps}, Bull. London Math. Soc.
10 (1978), 1-68.
\bibitem{4}
M. J. Gotay, J. Isenberg, J. E. Marsden, {\it Momentum Maps and the Hamiltonian
Structure of Classical Relativistic Fields}, http://xxx.lanl.gov/hep/9801019,
1998.
\bibitem{5}
N. Kamron et P. J. Olver, {\it Le Probl\'eme d'equivalence \`a une divergence
pr\'es dans le calcul des variations des int\'egrales multiples}, C. R. Acad.
Sci. Paris, t. 308, S\'erie I, p. 249-252, 1989.
\bibitem{6}
R. Miron, M. Anastasiei, {\it The Geometry of Lagrange Spaces: Theory and
Applications}, Kluwer Academic Publishers, 1994.
\bibitem{7}
R. Miron, M. S. Kirkovits, M. Anastasiei, {\it A Geometrical Model for Variational
Problems of Multiple Integrals}, Proc. of Conf. of Diff. Geom. and Appl. ,
June 26-July 3, 1988, Dubrovnik, Yugoslavia.
\bibitem{8}
M. Neagu, {\it Generalized Metrical Multi-Time Lagrangian Geometry of
Physical Fields}, Workshop on Diff. Geom. , Global Analysis, Lie Algebras, Aristotle University
of Thessaloniki, Greece, Aug. 27-Sept. 2, 2000; http://xxx.lanl.gov/math.DG/0011003, 2000.
\bibitem{9}
M. Neagu, {\it Harmonic Maps between Generalized Lagrange Spaces}, Southeast
Asian Bulletin of Mathematics, Springer-Verlag, 2000-2001, in press.
\bibitem{10}
M. Neagu, {\it Metrical Multi-Time Lagrangian Geometry of Physical Fields},
Workshop on Diff. Geom. , Global Analysis, Lie Algebras, Aristotle University
of Thessaloniki, Greece, Aug. 27-Sept. 2, 2000;
http://xxx.lanl.gov/math.DG/0009117, 2000.
\bibitem{11}
M. Neagu, {\it The Geometry of Relativistic Rheonomic Lagrange Spaces},
http://xxx.lanl.gov/math.DG/0010090, 2000.
\bibitem{12}
P. J. Olver, {\it Applications of Lie Groups to Differential Equations},
Springer-Verlag, 1986.
\bibitem{13}
D. Saunders, {\it The Geometry of Jet Bundle}, Cambridge University Press,
New York, London, 1989.
\bibitem{14}
Z. Shen, {\it Geometric Methods for Second Order Ordinary Differential
Equations}, preprint, 2000.
\bibitem{15}
C. Udri\c ste, {\it Nonclassical Lagrangian dynamics and potential maps}, Proc.
of the Conference on Mathematics in Honour of Prof. Radu Ro\c sca at the Occasion
of his Ninetieth Birthday, Katholieke University Brussel, Katholieke University
Leuven, Belgium, Dec. 11-16, 1999; http://xxx.lanl.gov/math.DS/0007060, 2000.
\bibitem{16}
C. Udri\c ste, {\it Solutions of DEs and PDEs as potential maps using first
order Lagrangians}, Centenial Vr\^anceanu, Romanian Academy, University of
Bucharest, June 30-July 4, 2000; http://xxx.lanl.gov/math.DS/0007061, 2000.
\bibitem{17}
C. Udri\c ste, M. Neagu, {\it Geometrical Interpretation of Solutions of
Certain PDEs}, Balkan Journal of Geometry and Its Applications, 4,1 (1999),
145-152.
\bibitem{18}
C. Udri\c ste, M. Neagu, {\it Metrical Multi-Time Lagrangian Interpretation
of PDEs}, 2000, to appear.
\end{thebibliography}
\end{document}